\documentclass{article}
\usepackage{bm,xspace,pstricks,amssymb,amsmath,graphicx}
\usepackage[dvips,colorlinks,bookmarksopen,bookmarksnumbered,citecolor=red,urlcolor=red]{hyperref}
\usepackage[right=3cm,left=3cm,top=3cm,bottom=3cm]{geometry}

\newcommand\BibTeX{{\rmfamily B\kern-.05em \textsc{i\kern-.025em b}\kern-.08em
T\kern-.1667em\lower.7ex\hbox{E}\kern-.125emX}}
\newcommand{\ie}{{\em i.e.\/}\xspace}
\newcommand{\eg}{{\em e.g.\/}\xspace}

\title{Bayesian models for cost-effectiveness analysis in the presence of structural zero costs}
\author{Gianluca Baio}
\date{\small \it Department of Statistical Science\\[-5pt] University College London (UK)}

\begin{document}
\setlength{\baselineskip}{18pt}
\maketitle

\begin{center}
{\fontsize{14}{16}\selectfont \textbf{Abstract}}
\vspace{20pt}

{\fontsize{9}{12}\selectfont
\begin{minipage}[c]{12cm}
Bayesian modelling for cost-effectiveness data has received much attention in both the health economics and the statistical literature in recent years. Cost-effectiveness data are characterised by a relatively complex structure of relationships linking the suitable measure of clinical benefit (\eg QALYs) and the associated costs. Simplifying assumptions, such as (bivariate) normality of the underlying distributions are usually not granted, particularly for the cost variable, which is characterised by markedly skewed distributions. In addition, individual-level datasets are often characterised by the presence of structural zeros in the cost variable. 

Hurdle models can be used to account for the presence of excess zeros in a distribution and have been applied in the context of cost data. We extend their application to cost-effectiveness data, defining a full Bayesian model which consists of a selection model for the subjects with null costs, a marginal model for the costs and a conditional model for the measure of effectiveness (conditionally on the observed costs). The model is presented using a working example to describe its main features.

\vspace{10pt}
\end{minipage}}
\end{center}
\vspace{10pt}

\section{Introduction}
Modelling for cost-effectiveness data has received much attention in both the health economics and the statistical literature in recent years \cite{Briggsetal:2006,WillanBriggs:2006}, increasingly often under a Bayesian statistical approach \cite{O'HaganStevens:2001,O'Haganetal:2001,Spiegelhalteretal:2004,Baio:2012}. From the statistical point of view, this is an interesting problem because of the generally complex structure of relationships linking a suitable measure of clinical benefit (\eg QALYs) and the associated costs. In addition, simplifying assumptions, such as (bivariate) normality of the underlying distributions are usually not granted, particularly for the cost variables.

In fact, costs are typically characterised by a markedly skewed distribution, which is generally due to the presence of a small proportion of individuals incurring large costs. To accommodate this feature, several models have been suggested and implemented. Among them, the most popular are probably represented by the log-Normal and Gamma distributions \cite{ThompsonNixon:2005,NixonThompson:2005}, which are well suited to describe right skewed data. 

However, in addition, individual-level datasets (such as those collected in clinical trials) are often characterised by the presence of \textit{structural} zeros in the cost variable: this amounts to observing a proportion of subjects for whom the observed cost is equal to zero. Under these circumstances, the use of log-Normal or Gamma models becomes impractical, since these distributions are defined for strictly positive arguments. A simple solution is to add a small constant $\epsilon$ to the entire set of observed values for the cost variable, thus artificially re-scaling it in the open interval $(0,\infty)$ \cite{Cooperetal:2007}. While very easy to implement, this strategy is potentially problematic, since the results are likely to be strongly affected by the actual choice of the scaling parameter $\epsilon$. In particular, there is no real guidance as to ``how small'' the value $\epsilon$ should be in order to minimise its influence on the economic results.

Alternatively, it is possible to use specific strategies to model data including structural zero costs that overcome this issue, \eg \textit{hurdle} models \cite{Ntzoufras:2009}; extensive treatment of this topic in the health economics literature is given in \cite{Cooperetal:2007,Cooperetal:2003,Mihaylovaetal:2011}, while applications include \cite{Toozeetal:2002,Harkanenetal:2013}. In a nutshell, the idea is to build a ``selection'' model that predicts the probability of a given individual being associated with a null cost; this is typically done using a logistic regression as a function of a set of relevant covariates. Then, for the individuals incurring a positive cost, a regression model is fitted to estimate the average cost, which effectively is a mixture of the two components.

With the notable exception of \cite{Lambertetal:2008} (who applied a bivariate Normal model to estimate survival and partially measured costs), hurdle models have been mainly used to either estimate the effect of relevant covariates or to predict future costs, without explicit reference to a measure of clinical benefit. The evaluation of the costs, however, is only one side of a comprehensive cost-effectiveness analysis, which needs to simultaneously account for the expected clinical benefits as well. As mentioned earlier, because costs and benefits are typically correlated, it is necessary to produce a multivariate model that can cater for this situation.

In this paper we aim at extending the two-parts model to produce a general framework able to account for: a) structural zero costs; and b) correlation between costs and clinical benefits. We take advantage of the flexibility of Bayesian models, which allow to specify several components that can then be linked to induce correlation among the different modules. We consider three components; the first one is a selection model that predicts the probability that each individual is associated with zero costs. The second module is a marginal model for the costs, which is expressed as a mixture of two components, depending on the observed value for the costs. Finally, the third module is a conditional model for the variable of effectiveness, conditionally on the observed value for the costs. 

The paper is structured as follows: first in section \ref{model} we set out our modelling framework. We then present the data and specific model used to analyse a case study in section \ref{example}, discussing the specific model in \ref{example_model} and the results in \ref{results}. Section \ref{conclusions} reviews the main conclusions.

\section{Modelling framework} \label{model}
Consider a dataset $\mathcal{D}$ including information on a set of $n$ individuals. This may arise in the case of a randomised clinical trial, or from observational data obtained from registries of clinical practice. We assume that, for each subject, $\mathcal{D}$ contains at least two variables $(e,c)$ measuring suitably defined clinical benefits and the associated costs. As we will show in the following, it is helpful to assume that the study also records some additional information at the individual level, \eg age, sex or potential co-morbidities. We also note that, even in the case of RCTs where these variables are not essential in the estimation of the treatment effects (by virtue of randomisation), they are usually measured and included in the final~dataset. 

For each intervention or treatment $t=0,\ldots,(T-1)$ under consideration, we can define a subset $\mathcal{D}_t$ with sample size $n_t$, so that 
$\mathcal{D}=\bigcup_t \mathcal{D}_t$ and $n=\sum_{t} n_t$. We partition the observed data as $\mathcal{D}_t = (\mathcal{D}_t^{\rm{null}} \cup \mathcal{D}_t^{\rm{pos}})$, where $ \mathcal{D}_t^{\rm{pos}}$ includes the $n_t^{\rm{pos}}$ individuals generating a positive cost. Consequently, $n_t^{\rm{null}}=n_t-n_t^{\rm{pos}}$ is the number of subjects with structural zero cost. Without loss of generality, we assume in the following that only two interventions are being considered: $t=0$ is some standard (\eg currently recommended or applied by the health care provider) and $t=1$ is a new intervention being suggested to potentially replace the standard.

\subsection{Selection model for $c=0$}\label{selectionModel}
We estimate the probability that each individual has a null observed cost, as a function of $J$ relevant covariates. For each subject in $i=1,\ldots,n_t$, we define an indicator $d_{it}$ taking value 1 if that individual is observed to have a null cost, and 0 otherwise. We model this variable as
\begin{eqnarray}
d_{it} & \sim & \mbox{Bernoulli}(\pi_{it}) \nonumber \\
\mbox{logit}(\pi_{it}) & = & \beta_{0t} + \sum_{j=1}^J \beta_{jt} Z^{t}_{ij}, \label{logistic}
\end{eqnarray}
where $\pi_{it}$ indicates the required probability. Both for computational and practical reasons (which we describe later), it generally helps to centre the covariates, \ie for each treatment group, instead of the originally observed covariate $\bm{X}^{t}_{j}$ we include in (\ref{logistic}) its centred version $Z^{t}_{ij} = X^{t}_{ij} - \mbox{E}[\bm{X}^{t}_{j}]$. Of course this construction implies that $\mbox{E}[\bm{Z}^t_j]=0$.

Within the Bayesian framework, the coefficients $\bm\beta_{t}=(\beta_{0t},\beta_{1t},\ldots,\beta_{Jt})$ are given suitable probability distributions. One simple choice is to use independent minimally informative Normal distributions (\ie centred on 0 and with large variance), but of course other choices are possible. This, however, does not normally impact to a (too) large extent on the results, especially in presence of (at least) moderately large datasets. 

Nevertheless, it is worth noting that, in cases where the number of subjects with structural zeros is very small, \textit{separation} (\ie the fact that a linear combination of the predictors is perfectly predictive of the outcome) is potentially an issue. A possible solution is to model the coefficients using Cauchy priors centred on 0 and with a small scale parameter \cite{Gelmanetal:2008}.

Under the assumptions specified above, the quantity
\[ p_{t} = \frac{\exp(\beta_{0t})}{1+\exp(\beta_{0t})} \]
represents the estimated overall probability of having a null cost for the ``average'' individual (\ie one with the values of the covariates set to 0, their mean). Sub-group analyses would be possible by selecting the combination of modalities for the covariates that define the required individual profile. Moreover, the model in (\ref{logistic}) can be extended to include individual structured (``random'') effects, for example in the case of clustering over time in repeated measurement data.

\subsection{Marginal model for the costs}\label{modelCost}
In the second module, we model the observed costs by specifying a single distributional form for the two components (subjects with null or positive costs). Nevertheless, this distribution is indexed by two different sets of parameters $\bm\theta_t=(\bm\theta_t^{\rm{null}},\bm\theta_t^{\rm{pos}})$, which depend on the value taken by $d_{it}$. In particular, we have that
\begin{eqnarray*} 
c_{it} \mid d_{it} & \sim & \left\{ \begin{array}{ll} p(c_{it} \mid d_{it}=1) = p(c_{it}\mid \bm\theta_t^{\rm{null}}), &\quad \mbox{for $i\in \mathcal{D}_t^{\rm{null}}$} \\ p(c_{it} \mid d_{it}=0) = p(c_{it}\mid \bm\theta_t^{\rm{pos}}), &\quad \mbox{for $i\in \mathcal{D}_t^{\rm{pos}}.$} \end{array} \right. 
\end{eqnarray*}
At this stage, we can choose any suitable distributional assumption for $p(c_{it}\mid d_{it})$.

\subsubsection{Gamma model}\label{gammaModel}
For example, we can model the costs for both components using a Gamma distribution
\[ c_{it}\mid d_{it} \sim \mbox{Gamma}\left(\eta_{t,d_{it}},\lambda_{t,d_{it}}\right), \]
where the nested index $d_{it}$ takes values 0,1 for patients with positive and null costs respectively. Thus, $\boldsymbol\theta_t^{\rm{pos}}=(\eta_{t0},\lambda_{t0})$ and $\boldsymbol\theta_t^{\rm{null}}=(\eta_{t1},\lambda_{t1})$, where, for $s,t=0,1$, $\eta_{ts}$ is the \textit{shape}, $\lambda_{ts}$ is the \textit{rate} of the Gamma density and, by definition,
\begin{eqnarray}
\psi_{ts} = \frac{\eta_{ts}}{\lambda_{ts}} \qquad \mbox{and} \qquad \zeta_{ts} = \sqrt{\frac{\eta_{ts}}{\lambda_{ts}^2}} \label{relGamma}
\end{eqnarray}  
are, respectively, the mean and the standard deviation of the cost distribution, on the natural scale.

It is possible to encode the required distributional assumptions by choosing suitable priors for $\boldsymbol\theta_t^{\rm{null}}$; in particular, we require that the distribution of costs must be identically 0 for all the patients with observed null cost. But if we choose $\eta_{t1}=w$ and $\lambda_{t1}=W$, with $0<w<<W$ (\eg $w=1$ and $W=10000$), then $\psi_{t1} \rightarrow 0$ and more importantly $\zeta_{t1} \rightarrow 0$ to an even faster rate of convergence. Thus, effectively, for the patients with observed null value, the cost is estimated to be identically~0 --- and as a matter of fact, this prior is so informative that no amount of evidence can modify it in the posterior. 

As for the patients with positive costs, we need to assume a non-degenerate prior on $\boldsymbol\theta_t^{\rm{pos}}$ to obtain a reasonable model. In this case, it is easier to encode any available information on the normal scale parameters $(\psi_{t0},\zeta_{t0})$, rather than on the shape and rate of the Gamma distribution, which are less straightforward to interpret and give a prior to. Just as an example, one may assume $\psi_{t0} \sim \mbox{Uniform}(0,H_\psi)$ and $\zeta_{t0} \sim \mbox{Uniform}(0,H_\zeta)$ for suitably selected values $H_\psi,H_\zeta$. Inverting the deterministic relationships in (\ref{relGamma}) it is easy to derive 
\[ \eta_{t0} = \psi_{t0} \lambda_{t0} \qquad \mbox{and} \qquad  \lambda_{t0} = \frac{\psi_{t0}}{\zeta_{t0}^2} \]
and thus the distributions selected for $(\psi_{t0},\zeta_{t0})$ automatically imply the priors for $(\eta_{t0},\lambda_{t0})$. Notice that these will in general not be vague at all, even in case the priors for $(\psi_{t0},\zeta_{t0})$ are chosen to be minimally informative, as in the example below --- in fact by assuming a flat prior on the natural scale parameter, we are implying some information on the orginal parameters of the assumed Gamma distribution. 

\subsubsection{Log-Normal model}\label{lognormal}
The above construction allows us to use a single distributional form to model the costs in both sub-groups and we can use the same rationale to encode different distributional assumptions. For example, we could use a log-Normal model to describe sampling variations in the observed costs. In this case, we have that 
\[ c_{it} \mid d_{it} \sim \mbox{log-Normal}\left(\eta_{t,d_{ti}},\lambda_{t,d_{it}}\right), \]
where $\eta_{ts}$ and $\lambda_{ts}$ are now the population average and standard deviation of the cost on the log scale, respectively. By the basic properties of the log-Normal distribution the mean and the standard deviation of the cost on the natural scale can be computed for each sub-group $s=0,1$ as
\begin{eqnarray} 
\psi_{ts} = \exp\left(\eta_{ts}+\frac{\lambda^2_{ts}}{2}\right) \qquad \mbox{and} \qquad \zeta_{ts} = \sqrt{\left(\exp(\lambda^2_{ts})-1\right)\exp\left(2\eta_{ts}+\lambda^2_{ts}\right)}. \label{relLogNorm}
\end{eqnarray}

Setting $\eta_{t1} = -W$ and $\lambda_{t1} = w$ for $0<w<<W$ (\eg $w=10$ and $W=1000$) implies that for the individuals in $\mathcal{D}_t^{\rm{null}}$ the cost \textit{on the natural scale} is effectively 0 with no substantial variability. Conversely, in a similar fashion to what shown above, we can define the non-degenerate prior for the individuals with positive costs on the natural scale, rather than on the original parameters on the log scale; for instance, one could again use Uniform priors for $(\psi_{t0},\zeta_{t0})$ and then invert the relationships in~(\ref{relLogNorm}) to obtain
\[ \eta_{t0} = \log(\psi_{t0}) - \frac{1}{2}\log\left[ 1+ \left(\frac{\zeta_{t0}}{\psi_{t0}}\right)^2 \right] \qquad \mbox{and} \qquad  \lambda_{t0} = \sqrt{\log\left[ 1+ \left(\frac{\zeta_{t0}}{\psi_{t0}}\right)^2 \right]} \]
and thus induce the priors for $(\eta_{t0},\lambda_{t0})$.

\subsubsection{Computation of the average cost}
Of course none of the distributional assumptions discussed above are essential and it is possible to express the available prior information in different ways and using other parametric models. Nevertheless, the general framework still applies and one can use a single distribution to represent the observed costs in both the components of the population, simply by cleverly modelling the parameters.

In addition, regardless on the underlying marginal model, once the two components of $\boldsymbol\psi_t=(\psi_{t0},\psi_{t1})$ have been estimated, it is then possible to derive the overall average cost in the population by computing the weighted average
\begin{eqnarray} 
\mu_{ct} & = & (1-p_t)\psi_{t0} + p_t\psi_{t1} \label{mixtureModel}\\
& = & (1-p_t)\psi_{t0} \nonumber ,
\end{eqnarray}
where the weights are given by the estimated probability associated with each of the two classes. In effect, the population average cost is obtained by down-weighting the estimated average for $\mathcal{D}_t^{\rm{pos}}$, to account for the presence of the structural zeros. The weights of the mixture are informed by the selection model.

\subsection{Conditional model for the measure of clinical benefit}
The final module consists in modelling the measure of clinical benefit $e$ so that correlation between the two dimensions of the health economic evaluation is accounted for. One possible way of doing so is to factorise the joint distribution $p(e,c\mid \boldsymbol\theta_t)$ in the product of a marginal and a conditional distribution. Intuitively, it is easier to think of this factorisation in terms of $p(e\mid \boldsymbol\theta_t)p(c\mid e,\boldsymbol\theta_t)$, \ie assuming that the observed costs somehow depend on the value taken by the measure of effectiveness. This construction has been used for example in \cite{O'HaganStevens:2001}.

However, because we are merely modelling a \textit{probabilistic} structure (\ie we are not claiming any \textit{causal} relationship), it is equally reasonable to factorise the joint distribution in terms of a marginal density for the costs and a conditional density for the benefits given the costs, \ie $p(e,c\mid \boldsymbol\theta_t)=p(c\mid \boldsymbol\theta_t)p(e\mid c,\boldsymbol\theta_t)$ --- in this sense we refer to the model of section \ref{modelCost} as \textit{marginal} and to the one in the current section as \textit{conditional}.

The distribution $p(e\mid c,\boldsymbol\theta_t)$ is chosen according to the nature of the effectiveness variable. For example, if $e$ were expressed in terms of QALYs over a long period of time, it should be a continuous density defined in $\mathbb{R}^+$. But, as discussed in \cite{ThompsonNixon:2005}, whatever this choice one can always represent its mean (\ie the conditional average effectiveness, given the costs) through a regression model
\begin{eqnarray} 
g(\phi_{it}) = \xi_t + \gamma_t(c_{it}-\mu_{ct}) \label{regr_E_C}
\end{eqnarray}
defined in terms of a suitable link function $g(\cdot)$. The form of the link function obviously depends on the scale in which $\phi_{it}$ is defined, \eg if $\phi_{it}$ were modelled on the natural scale of $e$, then $g(\cdot)$ would be the identity function. 

In equation (\ref{regr_E_C}), the coefficient $\mu_{ct}$ is the population average cost obtained in the mixture model of (\ref{mixtureModel}), while the coefficients $\xi_t$ and $\gamma_t$ represent respectively the population (marginal) average effectiveness, and the level of correlation between effectiveness and costs. Notice that these are quantified on the scale defined by the link function. Thus, in order to estimate the marginal average effectiveness on the natural scale, it is necessary to compute the inverse transformation $\mu_{et} = g^{-1}(\xi_t)$. To complete the full Bayesian model, the parameters $(\xi_t,\gamma_t)$ as well as any other nuisance parameter characterising $p(e\mid c,\boldsymbol\theta_t)$ are given appropriate prior distributions. 

\begin{figure}[!h]
\centering
\includegraphics*{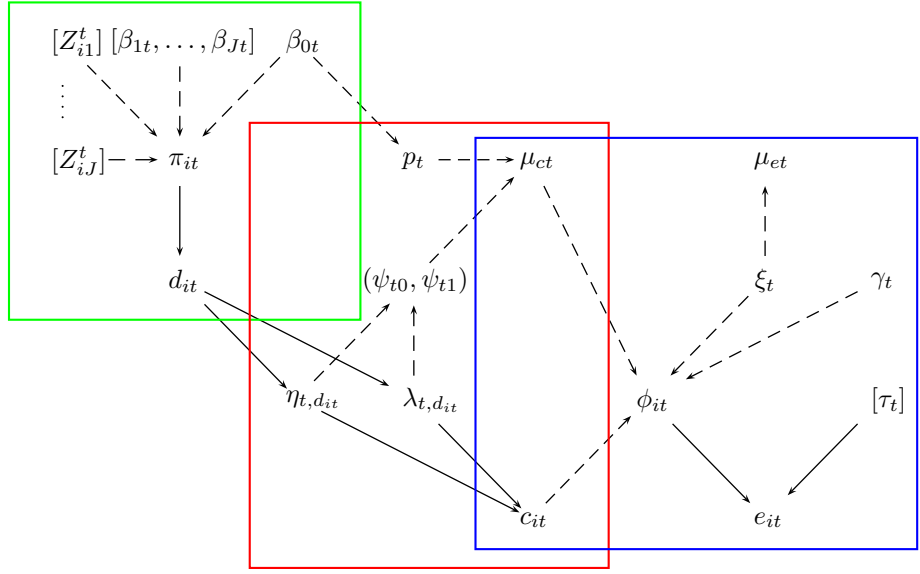}
\caption{A graphical representation of the full Bayesian model accounting for: a) the selection model for $c=0$; b) the marginal model for the costs; and c) the conditional model for the clinical benefit (given the observed costs). Dashed connections indicate logical relationships, while solid connections indicate probabilistic relationships. Nodes enclosed in brackets may be not be used: for example, the covariates $\bm{Z}^t_1,\ldots,\bm{Z}^t_J$ may not be observed and hence the coefficients $\beta_{1t},\ldots,\beta_{Jt}$ are not included in the selection model; similarly, the parameter $\tau_t$ may not be needed in the conditional model for $e$ (\eg if a Bernoulli distribution is considered, only $\phi_{it}$ is necessary)}\label{GraphModel}
\end{figure}

Figure \ref{GraphModel} shows a graphical representation of the general model structure, highlighting the links among the three modules. Dashed connections indicate logical relationships among nodes (variables), while solid connections represent probabilistic relationships or dependence. The three modules are connected by means of these relationships.

\subsection{Economic evaluation}
Once the model is fitted to the observed data $\mathcal{D}$, it is possible to directly use the posterior distributions for $(\mu_{et},\mu_{ct})$ to perform the health economic evaluation. For example, we can construct suitable health economic summaries, such as the \textit{increment in mean effectiveness} $\Delta_e := \mu_{e1} - \mu_{e0}$ and the \textit{increment in mean cost} $\Delta_c := \mu_{c1} - \mu_{c0}$.
 
After having obtained the required posteriors, for instance using an MCMC procedure, one can post-process the output (\eg using the \texttt{R} package \texttt{BCEA} \cite{Baio:2012,Baio:2013a}) and perform the economic analysis. This includes constructing the cost-effectiveness plane, which describes the joint distribution of $(\Delta_e,\Delta_c)$, and the expected incremental benefit $\mbox{EIB}=k\mbox{E}[\Delta_e]-\mbox{E}[\Delta_c]$, which is used to perform the decision analysis upon deterministically varying the willingness-to-pay threshold $k$\footnote{The willingness-to-pay $k$ is used to put the cost and effectiveness differentials on the same scale and it represents the cost that the decision maker is willing to pay to increment the effectiveness measure by one unit. If $\mbox{EIB}>0$ then, for a given $k$, $t=1$ is more cost-effective than $t=0$. More details are presented in \cite{Briggsetal:2006,WillanBriggs:2006,O'HaganStevens:2001,O'Haganetal:2001,Spiegelhalteretal:2004,Baio:2012}.}. 

Moreover, it is helpful to conduct a probabilistic sensitivity analysis (PSA), \eg in terms of the cost-effectiveness acceptability curve $\Pr(k\Delta_e-\Delta_c>0)$ and the analysis of the expected value of information \cite{Spiegelhalteretal:2004,Baio:2012}, in order to assess the impact of the parameters uncertanty on the decision process. Finally, as different distributional assumptions for the costs may also have an impact on the decision-making, it is generally advisable to perform a structural sensitivity analysis \cite{Jacksonetal:2011}. Our framework allows these to be done in a straightforward way. 

\subsection{Sensitivity to the parameters specification for $(\psi_{t1},\zeta_{t1})$}
The choice of the values $w$ and $W$ in the models for $(\psi_{t1},\zeta_{t1})$ described above is potentially a delicate issue, as the results may be sensitive to their specification. In fact, the estimation for the main parameters $(\mu_{et},\mu_{ct})$ is not really affected by this choice, provided that the encoded relationship between them really induces $\psi_{t1},\zeta_{t1}\rightarrow 0$.

It is worth noticing that, as is reasonable, different values for $w$ and $W$ do have an impact on measures of model fit, such as the Deviance Inflation Criterion (DIC) \cite{Spiegelhalteretal:2002}. This is essentially due to the fact that the population is really made by two groups, one of which shows costs that are identically null. Thus, the faster the rate of convergence to 0 for $\zeta_{t1}$, the better the fit to the observed data and therefore the smaller the resulting DIC. But the implications in terms of the resulting estimated values (and hence the resulting economic model) are immaterial.

\section{Example: acupuncture trial}\label{example}
We consider the RCT conducted in the UK primary care setting to evaluate the cost-effectiveness of acupuncture in the management of chronic headache, originally presented in \cite{Wonderlingetal:2004}. The trial recruited a total of 401 patients aged 18-65 who reported an average of at least two headaches per month, from general practices in England and Wales. The participants were randomly allocated to either usual care (which we indicate with $t=0$), or in addition up to 12 acupuncture treatments over three months from appropriately trained physiotherapists (active intervention, $t=1$).

The measure of effectiveness used in the study is the total QALYs gained, obtained using a specific algorithm based on the SF-6D questionnaire \cite{Brazieretal:2002}. Only 255 patients ($n_0=119$ in the control and $n_1=136$ in the active treatment group) had valid data for the SF-6D questionnaire and thus the economic evaluation is performed on this sub-sample. Notice that because the time horizon considered is one year and the SF-6D utility measure is defined in the interval $[0;1]$, the resulting QALYs for each individual are also constrained in this interval. The overall cost was calculated for each patient by adding up the several resources, including non-prescription drugs and private health care visit, visits to practitioners of complementary or alternative medicine and drug prescriptions. 

\subsection{Model specification}\label{example_model}
We follow the framework of section \ref{model} to model the indicator of structural zero cost; however, because no additional covariate at the individual level is available, we only use the intercept $\beta_{0t}$ to predict the probability of zero cost. To make a parallel with the missing data literature, we are then effectively assuming a mechanism of ``zero completely at random'' (ZCAR), in which the chance of observing an individual associated with zero costs does not depend on any other variable\footnote{In fact, the correspondence between missing data (MD) and structural zeros (SZ) is not perfect. In MD, the assumption of missing completely at random (MCAR) implies that the ``missingness'' and the ``outcome'' modules are completely separated (and thus the chance of observing a missing value is assumed to be independent on any other variable, including the outcome). This implies that under MCAR there is no need to include the missingness module in the analysis. On the other hand, in SZ the two modules are always linked, because the distribution of the outcome depends on the zero indicator. Nevertheless, under ZCAR, we are assuming the absence of other (observed or unobserved) factors that can influence the chance of observing a zero.}. Of course, this may or may not be appropriate and typically a less restrictive assumption of ``zero at random'' (ZAR), under which the zero pattern depends only on observable covariates (as in Figure \ref{GraphModel}), will be more tenable. In order to make the estimation more stable, we model $\beta_{0t}\sim\mbox{Cauchy}(0,2.5)$.

As for the costs, we use both a Gamma and a log-Normal specification, as described in sections \ref{gammaModel} and \ref{lognormal}. In the Gamma model, we set $w=1$ and $W=10000$, while in the log-Normal model, we specify $w=1$ and $W=50$. In both cases, we set $H_\psi=1\,000$ and $H_\zeta=300$; these values encode the knowledge that the intervention is not particularly expensive, and perform a sensitivity analysis to the choice of the values for the parameters $(w,W)$.

Finally, we model the effectiveness measure using a Beta regression, which in line with \cite{Figueroaetal:2013} we specify as follows
\begin{eqnarray}
e_{it} \mid c_{it} & \sim & \mbox{Beta}\left(\phi_{it}\tau_t,(1-\phi_{it})\tau_t\right) \label{effect} \\
\mbox{logit}(\phi_{it}) & = & \xi_t + \gamma_t(c_{it}-\mu_{ct}) \nonumber \\
\xi_t,\gamma_t,\log(\tau_t) & \stackrel{iid}{\sim} & \mbox{Normal}(0,10\,000), \nonumber 
\end{eqnarray}
In equation (\ref{effect}), the parameter $\phi_{it}$ represents the conditional subject-specific average QALYs, while the parameter $\tau_t$ is the conditional precision (inverse variance), which we assume constant across the subjects within each treatment arm. The actual measure of effectiveness (\ie the marginal population average QALYs under either treatment) can be then retrieved on the correct scale by applying the inverse logit transformation
\[ \mu_{et} = g^{-1}(\xi_t) = \frac{\exp(\xi_t)}{1+\exp(\xi_t)}. \]

\subsection{Results}\label{results}
We fitted the models of section \ref{example_model} using the \texttt{R} package \texttt{BCEs0} \cite{Baio:2013b}, which implements the general framework described in section \ref{model} under a set of possible distributional assumptions. In \texttt{BCEs0}, the user needs to: \textit{i}) specify a data list including the observed values for $(e,c)$ under the two treatment options and possibly the matrices including the values for the covariates $\bm{X}^t_1,\ldots,\bm{X}^t_J$, and the fixed parameters $H_\psi$ and $H_\zeta$; \textit{ii}) select a distribution for the costs (implemented choices are Gamma, log-Normal and Normal); and \textit{iii}) select a distribution for the measure of effectiveness (Beta, Bernoulli, Gamma and Normal are currently implemented). \texttt{BCEs0} will then write the \texttt{JAGS} \cite{jags} model for the selected specification to a text file, call the library \texttt{R2jags} \cite{R2jags} (which connects \texttt{JAGS} to \texttt{R} in background) and perform the MCMC analysis. The resulting simulations from the posterior distributions are saved to the \texttt{R} workspace and can be used for the health economic evaluation. Since the model file is saved in the working directory, the user can use it as a template and modify.

We ran 10\,000 iterations, using a burn-in of 5\,000 and retaining one iteration every 10, which resulted in a sample of 1\,000 iterations which we used to produce the posterior analysis. For each variable in the model, convergence of the MCMC sampler was assessed by the analysis of the potential scale reduction \cite{GelmanRubin:1992}, as well as the effective sample size. 

Table \ref{table} presents summary statistics from the posterior distributions of the main parameters in the model, for both specifications of the cost variable. In both models, treatment $t=1$ is associated with both higher costs and higher QALYs, on average. In particularly, the average costs is substantially larger for this arm of the trial. As is possible to see, for both treatments there is a significant difference between $\mu_{ct}$, the overall average cost and $\psi_{t0}$, the average cost for the subjects in~$\mathcal{D}_t^{\rm{pos}}$. The log-Normal model produces estimations of the costs that are slightly lower than those produced by the Gamma specification. 

\begin{table}[!h]
\caption{Posterior summaries for selected parameters for the Gamma/Beta and log-Normal/Beta models}\label{table}
\centering
\begin{tabular}{crrrrrrrrrrr}
\hline
& \multicolumn{4}{c}{Gamma/Beta model}	&	& \multicolumn{4}{c}{log-Normal/Beta model} \\
\hline
Parameter	&Mean	&SD	& \multicolumn{2}{c}{95\% interval}	&	& Mean	&SD	& \multicolumn{2}{c}{95\% interval}\\
\hline
$p_0$	&0.039	&0.018	&0.012	&0.080 &	&0.038	&0.019	&0.011	&0.080\\
$\psi_{00}$	&226.958	&20.885	&187.075	&267.431 &	&179.934	&11.916	&155.702	&201.788\\
$\mu_{c0}$	&218.150	&20.520	&176.721	&259.649 &	&173.071	&11.830	&148.505	&195.024\\
$\mu_{e0}$	&0.710	&0.011	&0.688	&0.731 &	&0.712	&0.010	&0.691	&0.731\\[5pt]
$p_1$	&0.011	&0.009	&0.001	&0.036 & &0.011	&0.009	&0.001	&0.032\\
$\psi_{10}$	&408.171	&21.213	&365.324	&449.598 &	&393.165	&18.378	&357.348	&428.717\\
$\mu_{c1}$	&403.823	&21.411	&361.115	&445.206 &	&378.179	&18.972	&341.202	&415.554\\
$\mu_{e1}$	&0.729	&0.011	&0.708	&0.750 &	&0.731	&0.011	&0.710	&0.751\\
\hline
\end{tabular}
\end{table}

Of course, in a real-life, comprehensive analysis, sensitivity to the choice of the distribution for the parameters of each model should be explored extensively. Incidentally, the values for the willingness to pay parameter $k$ beyond which $\mbox{EIB}>0$ for the two different specifications of the model are quite similar (\pounds 9\,754 and \pounds 10\,713 for the Gamma/Beta and the log-Normal/Beta models, respectively) and so are the results of PSA. Thus, the practical implications of the resulting differences would not be particularly relevant in this case, as the active intervention would be considered as cost-effective and sustainable for the NHS under both specifications.

\begin{figure}
\centering
\includegraphics[scale=.55]{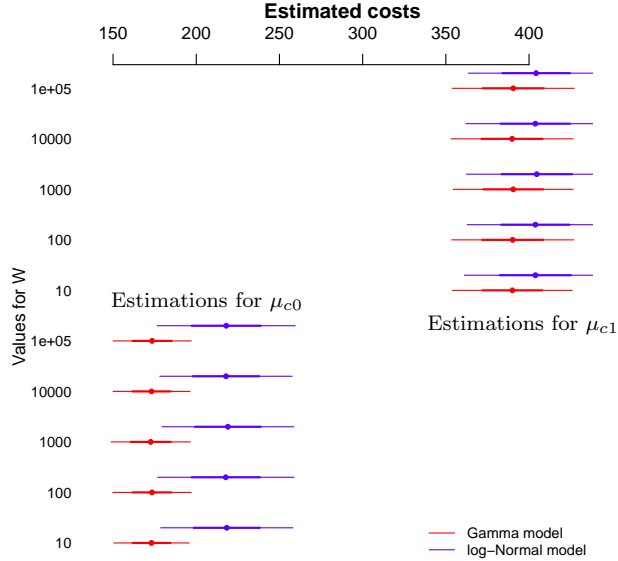}
\rput(-5.8,3.8){\fontsize{8}{8}\selectfont Estimations for $\mu_{c0}$}
\rput(-1.6,3.5){\fontsize{8}{8}\selectfont Estimations for $\mu_{c1}$}
\caption{Sensitivity analysis for the choice of the parameter $W$. The dots represent the posterior means for the estimated costs $\mu_{c0}$ and $\mu_{c1}$ (on the left and the right of the panel, respectively). The light dark and lines indicate the 50\% and 95\% credible intervals, respectively. Dots and lines in red indicate the Gamma model, while those in blue indicate the log-Normal model. Within models, in all cases the results are substantially identical and do not depend at all on the selected value of $W$}\label{sensW}
\end{figure}

We have run the models using different values for the parameters $w,W$, to assess the impact on the actual cost estimation. As an example, Figure \ref{sensW} shows the results of the sensitivity analysis on $W$, holding $w$ fixed to 1 for all cases; in particular, we have selected values of $W=(10,100,1\,000,10\,000,100\,000)$. In the graph, we report the posterior mean and both a 50\% and a 95\% posterior credible interval for the average costs (the dark and light lines, respectively). The results for $t=0$ are depicted on the left side, while those for $t=1$ are on the right side. As is possible to see, the point as well as the interval estimate of the average costs are effectively unchanged in all the cases. 

It is interesting to notice that model fitting (as measured in terms of the DIC) varies differently for the Gamma and the log-Normal model. In the former, the DIC becomes smaller when $W$ increases and, although there are no practical differences in terms of the estimation for all the parameters, the best-fitting model is the one with $W=100\,000$, although the results are hardly different for all the parameters. 

Conversely, for the log-Normal specification, the best-fitting model is the one with $W=10$; however, convergence is not really achieved in this case, particularly for the coefficient $\beta_{01}$, which estimates the chance of zero cost for $t=1$. In the absence of additional information, the model effectively computes a marginal average probability based only of the observed data (and no other covariate). But since the observed proportion is very close to 0, the estimations are unstable when $W$ is not large enough. We reiterate that this has a meaningful implication on the model \textit{convergence}, but not on the estimated costs and effectiveness, and hence on the health economic evaluation. A value of $W=50$ (which is used by default in \texttt{BCEs0}) is sufficient to provide convergence and model fit.

\section{Discussion}\label{conclusions}
In this paper we have defined and discussed a general framework to handle cost-effectiveness analysis using individual level data (\eg from a RCT) in the presence of structural zeros in the cost variable. This is a challenging situation because cost-effectiveness models are characterised by a relatively complex structure which require the formal inclusion of correlation between the outcomes. In addition, because of asymmetry in the cost distributions, we also need to model them using suitable formulations.

The framework developed in section \ref{model} uses of a flexible structure and allows the two components of the cost distributions to be modelled using a single specification. The parameters of the cost distributions are defined differently in the two components of the mixture; for the subjects with observed null costs, the specification implies that the final estimation is identically 0. Consequently, the final estimation of the overall population mean costs is a weighted average of the two components. The correlation between costs and clinical benefits is ensured by the model strucure.

The choice of the selection model for $c>0$ is of course crucial. In the example of section \ref{example}, we have assumed that no unobserved factors impact on the probability of an individual generating null costs. This is in general a very restrictive assumption and in fact it is likely that, even in a randomisation context, the probability $\pi_{it}$ does depend on some covariates and, possibly, is also affected by unobserved confounders. In general, this complicates the situation, but the generalisability of the framework is not compromised, as the model in (\ref{logistic}) can be extended to deal with more complex situations.

The \texttt{R} package \texttt{BCEs0} can be used, at least as a first approximation, to build a model consistent with the general framework. The choice of possible distributions for $(e,c)$ is limited to what we consider to be the most likely situations. However, a translation into the \texttt{JAGS} language (which is effectively identical when applied to other software such as \texttt{OpenBUGS}) is automatically generated. Thus, the user can easily modify the ``template'' model file to cater for their specific needs (\eg adding a different distribution for $e\mid c$, including structured effects in the selection model, or modifying the priors).

Unlike simpler but less efficient solutions to the problem of structural zeros in cost-effectiveness analysis, the structure of section \ref{model} is quite robust to the choice of the relevant parameters. In particular, it is pretty easy to ensure that the null component of the mixture for the cost distributions is indeed identically 0. As showed in the acupuncture trial example, the choice of the parameters $w,W$ is effectively irrelevant, provided that the mean and variances for the distribution under $\bm\theta_t^{null}$ tend to zero sufficiently quickly.

\section*{Acknowledgments}
The acupuncture trial data are kindly provided by David Wonderling, Richard Nixon and Richard Grieve.

\end{document}